\documentclass{amsart}
\usepackage{pstricks}

\newtheorem{theorem}{Theorem}
\title{A simple bijection between permutation tableaux and permutations}

\author{Sylvie Corteel}

\address{CNRS LRI, Universit\'e Paris-Sud, B\^atiment 490, 91405 Orsay
Cedex France} 
\email{Sylvie.Corteel@lri.fr}

\begin{document}

\maketitle 

Permutation tableaux are new objects that come
from the enumeration of  the totally positive Grassmannian cells
\cite{Postnikov, Williams}. Surprisingly  they are also connected
to a statistical physics  model called 
the Partially ASymmetric Exclusion Process \cite{Corteel,CW,CW1}.


Our main interest here is that these tableaux are in bijection
with permutations. To our knowledge only one bijection
between permutations and permutation tableaux is known and
appeared in \cite{SW}, although several bijections between
a subclass called Catalan tableaux and diverse objects counted
by Catalan numbers are known \cite{V}. This bijection between 
permutations and permutation tableaux is quite complicated;
but a lot of statistics of the permutation (weak exceedances, crossings \cite{Corteel},
alignments \cite{Williams}\ldots) can be read from the tableau. In particular
the number of weak exceedances of the permutation corresponds to the
number of rows of the tableau. See \cite{SW}
for many more details.

Our goal here is therefore to present a simple bijection
between permutations and permutation tableaux that
maps the descents of the permutation to the columns of the 
tableau. Our result is 
\begin{theorem}
There exists a bijection between permutations of $\{1,\ldots ,n\}$ with $k$ descents
and permutation tableaux
of length $n$ with $k$ columns.
\end{theorem}

We first start with a few definitions. 
A descent of a permutation $\sigma$ of a set $S$ is an entry $i\in S$ such
that if $\sigma(j)=i$ then $\sigma(j)>\sigma(j+1)$.
For example, if $\sigma=(2,4,8,5,1,6,3,7)$ and $S=\{1,2,3,4,5,6,7,8\}$, then 5,6 and 8 
are the descents of $\sigma$.

As in \cite{SW}, a {\em permutation tableau} $T$ is a shape (the
Ferrers diagram of a  partition) 
together with a filling of the cells with $0$'s and
$1$'s such that the following properties hold:
\begin{enumerate}
\item Each column contains at least one $1$.
\item There is no $0$ which has a $1$ above it in the same column
{\em and} a $1$ to its left in the same row.
\end{enumerate}
An example of a permutation tableau is given in Figure \ref{fig}.
Different statistics on permutation tableaux were defined in \cite{CW1,SW}.
We list a few here. The {\em length} of a tableau is the half perimeter of its shape.
A zero in a permutation tableau is  {\it restricted} if 
there is a one above it in the same column. A row is 
{\it unrestricted}
if it does not contain a restricted entry.  A restricted zero is a {\em rightmost}
restricted zero if it is restricted and it has no restricted zero to its right in the same row.
\begin{figure}[h]
 \centering
        \psset{unit=0.4cm}
\begin{pspicture}(0,0)(3,5)        
\psframe[dimen=middle](0,0)(1,1)
\rput(1.2,0.5){\tiny 7}
\rput(0.5,-0.2){\tiny 8}
\psframe[dimen=middle](0,1)(1,2)
\psframe[dimen=middle](1,1)(2,2)
\psframe[dimen=middle](2,1)(3,2)
\rput(3.2,1.5){\tiny 4}
\rput(1.5,0.8){\tiny 6}
\rput(2.5,0.8){\tiny 5}
\psframe[dimen=middle](0,2)(1,3)
\psframe[dimen=middle](1,2)(2,3)
\psframe[dimen=middle](2,2)(3,3)
\rput(3.2,2.5){\tiny 3}
\rput(3.2,3.5){\tiny 2}
\rput(3.2,4.5){\tiny 1}

\psframe[dimen=middle](0,3)(1,4)
\psframe[dimen=middle](1,3)(2,4)
\psframe[dimen=middle](2,3)(3,4)

\psframe[dimen=middle](0,4)(1,5)
\psframe[dimen=middle](1,4)(2,5)
\psframe[dimen=middle](2,4)(3,5)

\rput(0.5,0.5){\small 1}
\rput(0.5,1.5){\small 0}
\rput(1.5,1.5){\small 1}
\rput(2.5,1.5){\small 1}
\rput(0.5,2.5){\small 1}
\rput(1.5,2.5){\small 1}
\rput(2.5,2.5){\small 1}
\rput(0.5,3.5){\small 0}
\rput(1.5,3.5){\small 0}
\rput(2.5,3.5){\small 1}
\rput(0.5,4.5){\small 1}
\rput(1.5,4.5){\small 0}
\rput(2.5,4.5){\small 1}
\end{pspicture}
        \caption{Example of a permutation tableau}
\label{fig}
\end{figure}

We label the boundary 
of the shape of the tableau from 1 to its length, going from top-right
to bottom-left. See Figure \ref{fig}.
This labels the rows and the columns :
a South  (resp. West) step labelled by $i$ gives the label $i$ to this row (resp. column).
The cell $(i,j)$ of the tableau  corresponds the cell that is in the row labelled by $i$ and the column
labelled by $j$. 

In Figure \ref{fig}, a permutation tableau of shape $(3,3,3,3,1)$ and length 8 is given. 
The rows 1,3 and 7 are unrestricted and the rows 2 and 4 are restricted. 
The rightmost restricted zeros are in cells $(2,8)$ and $(4,8)$.

We are now ready to present the bijection.
Let $\sigma$ be a permutation of $\{1,\ldots ,n\}$ and let $T$ be its image.
We first draw the shape of the tableau $T$ . For $i$ from 1 to $n$, we draw
a West step if $i$ is a descent and a South step otherwise. An example for
$\sigma=(2,4,8,5,1,6,3,7)$ is given in Figure \ref{fig}, as 5,6 and 8 are the descents
of $\sigma$.

Now let us fill the cells of the tableau $T$. Let $(i,j)$ be the  Eastmost and Southmost cell
of the tableau $T$ that is not filled:
\begin{itemize}
\item if $i$ and $j$ are not adjacent in the permutation $\sigma$ then fill the cell $(i,j)$ with a one.
\item otherwise 
\begin{itemize}
\item if $i$ is before $j$, then fill all the empty cells of row $i$ with zeros and delete
$i$ from the permutation $\sigma$.
\item otherwise fill cell $(i,j)$ with a one and all the empty cells of column $j$ with zeros and delete
$j$ from the permutation $\sigma$.
\end{itemize}
\end{itemize}
At the end of this process, $T$ is filled and $\sigma$
is the list of the labels of the unrestricted rows of $T$ in increasing order. 
It is easy to check that $T$ is a permutation tableau and that $i$ is a 
descent in $\sigma$ if and only if there is a column labelled by $i$ in $T$.\\

\noindent {\bf Example 1.}  We start with $\sigma=(2,4,8,5,1,6,3,7)$, and we draw the shape of $T$ 
(see Figure \ref{fig}). We first fill cell $(4,5)$ and fill it with a 1, as
4 and 5 are not adjacent. We also fill cells $(3,5)$ and $(2,5)$ with ones. We fill (1,5) with a one
and delete 5 from the permutation. The permutation is now $\sigma=(2,4,8,1,6,3,7)$.
We fill cell $(4,6)$ with a one and then cell $(3,6)$ with a one and all the cells above in this column
with a zero and delete 6 from the permutation. The permutation is now $\sigma=(2,4,8,1,3,7)$. 
Then cell $(7,8)$
gets a one, cell $(4,8)$ gets a zero and 4 is deleted. 
The permutation is now $\sigma=(2,8,1,3,7)$. Finally cell $(3,8)$ gets a one, cell $(2,8)$ 
gets a zero and $2$ is deleted. 
The permutation is now $(8,1,3,7)$. Cell $(1,8)$ gets a one and all the cells above get a zero and
8 is deleted. The permutation is finally $(1,3,7)$. The result is given in Figure \ref{fig}.

\noindent {\bf Example 2.}
Starting with the permutation $(8,5,4,7,2,3,1,6)$, we first draw the shape of the tableau $T$
(Figure \ref{fig1}).
We first fill cell (2,3) and all the cells to its left with zeros and delete 2 from the permutation.
Then we add a one in cell (1,3) and delete 3 from the permutation. Then we fill cell (4,5) with  a one 
and all the cells above with zeros and delete 5 from the permutation. Then we fill cell (6,7) with a one
and fill cell (4,6) and all the cells to its left with zeros and delete 4 from the permutation. Then 
we fill cell (1,7) with a one and delete 7 from the permutation. Finally we fill cell (6,8) and cell (1,8)
with a one and delete 8 from the permutation. See Figure \ref{fig1}.

\begin{figure}[h]
 \centering
        \psset{unit=0.4cm}
\begin{pspicture}(0,0)(3,5)        
\psframe[dimen=middle](0,0)(1,1)
\psframe[dimen=middle](1,0)(2,1)
\rput(2.2,0.5){\tiny 6}
\rput(0.5,-0.2){\tiny 8}
\rput(1.5,-0.2){\tiny 7}
\psframe[dimen=middle](0,1)(1,2)
\psframe[dimen=middle](1,1)(2,2)
\psframe[dimen=middle](2,1)(3,2)
\rput(3.2,1.5){\tiny 4}
\rput(2.5,0.8){\tiny 5}
\psframe[dimen=middle](0,2)(1,3)
\psframe[dimen=middle](1,2)(2,3)
\psframe[dimen=middle](2,2)(3,3)
\psframe[dimen=middle](3,2)(4,3)

\rput(4.2,2.5){\tiny 2}
\rput(3.5,1.8){\tiny 3}
\rput(4.2,3.5){\tiny 1}

\psframe[dimen=middle](0,3)(1,4)
\psframe[dimen=middle](1,3)(2,4)
\psframe[dimen=middle](2,3)(3,4)
\psframe[dimen=middle](3,3)(4,4)

\rput(0.5,0.5){\small 1}
\rput(1.5,0.5){\small 1}

\rput(0.5,1.5){\small 0}
\rput(1.5,1.5){\small 0}
\rput(2.5,1.5){\small 1}

\rput(0.5,2.5){\small 0}
\rput(1.5,2.5){\small 0}
\rput(2.5,2.5){\small 0}
\rput(3.5,2.5){\small 0}

\rput(0.5,3.5){\small 1}
\rput(1.5,3.5){\small 1}
\rput(2.5,3.5){\small 0}
\rput(3.5,3.5){\small 1}
\end{pspicture}
        \caption{Image of the permutation $(8,5,4,7,2,3,1,6)$.}
\label{fig1}
\end{figure}

The reverse is as easy to define.  We start with the tableau $T$.
Then we initialize the permutation $\sigma$ to the list of the labels of the unrestricted rows in increasing order. Now for each column, starting from the left proceeding to the right, if the column is 
labelled by $j$
and if $(i,j)$ is the topmost one of the column then we add $j$ to the left of $i$ in the permutation
$\sigma$.
Moreover if column $j$ contains rightmost restricted zeros in rows $i_1,\ldots ,i_k$ then
we add $i_1,\ldots ,i_k$ in increasing order to the left of $j$ in the permutation $\sigma$.
It is easy to see that this is the reverse mapping.\\

\noindent{\bf Example 1.} We start with the tableau in Figure \ref{fig}. The unrestricted rows are rows
1,3 and 7.
The rightmost restricted zeros are in cells $(2,8)$ and $(4,8)$.
We start with the permutation $(1,3,7)$,
We add 8 to the left of 1 and add 2 and 4 to the left of 8. We get $(2,4,8,1,3,7)$. We add 6
to the left of 3 and get $(2,4,8,1,6,3,7)$. Finally we add 5 to the left of $1$. The permutation is $(2,4,8,5,1,6,3,7)$.

\noindent{\bf Example 2.} We start with the tableau in Figure \ref{fig1}. The unrestricted rows are rows 1 
and 6. The rightmost restricted zeros are in cells $(4,7)$ and $(2,3)$. We start with the permutation
$(1,6)$. We add 8 to the left of 1 and get $(8,1,6)$. We add 7 to the left of 1 and 4 to the left of 7
and get $(8,4,7,1,6)$. We then add 5 to the left of 4 and get $(8,5,4,7,1,6)$. Finally we add 3 to
the left of 1 and 2 to the left of 3. The result is $(8,5,4,7,2,3,1,6)$.

\end{document}